\documentclass[12pt]{article}

\usepackage{epsfig}
\usepackage{graphicx}
\usepackage{amsmath,amsthm}
\usepackage{mathbbol}
\usepackage{amssymb}
\usepackage{latexsym}
\usepackage[round]{natbib}
\usepackage{float}
\usepackage{color}
\usepackage{soul}
\usepackage{setspace}
\newcommand{\dd}{{\, \mathrm{d}}}
\newcommand{\bX}{\mathbf{X}}
\newcommand{\bY}{\mathbf{Y}}
\newcommand{\bZ}{\mathbf{Z}}
\newcommand{\bx}{\mathbf{x}}
\newcommand{\by}{\mathbf{y}}
\newcommand{\bfx}{\mathbf{x}}
\newcommand{\xn}{{\{x_1,\ldots,x_n\}}}
\newcommand{\xnl}{{x_1,\ldots,x_n}}
\newcommand{\bXred}{\mathbf{X}^!_\xnl}
\newcommand{\ta}{\theta}
\newcommand{\R}{\mathbb R}
\newcommand{\E}{\mathrm E}
\newcommand{\PP}{\mathrm P}
\newcommand{\CC}{\mathrm C}

\newcommand{\ld}{\lambda}

\usepackage{authblk}

\begin{document}
 
\title{A tutorial on Palm distributions for spatial point processes}

\author[1]{Jean-Fran{\c{c}}ois Coeurjolly}
\affil[1]{Laboratory Jean Kuntzmann, Statistics Department, Grenoble Alpes University,
51 Rue des Math{\' e}matiques,
Campus de Saint Martin d'H{\` e}res,
BP 53 - 38041 Grenoble cedex 09, France, email: \texttt{Jean-Francois.Coeurjolly@univ-grenoble-alpes.fr}.} 
\author[2]{Jesper M{\o}ller}
\author[2]{Rasmus Waagepetersen}  
\affil[1]{Department of Mathematical Sciences,
 Aalborg University, Fredrik Bajersvej 7E, DK-9220 Aalborg, email:
\texttt{jm@math.aau.dk, rw@math.aau.dk, }}

\maketitle

\begin{abstract}
This tutorial provides an  introduction to Palm distributions
  for spatial point processes. Initially, in the context of
  finite point processes, we give an 
explicit definition of Palm distributions in terms of their density
functions. Then we review Palm distributions in the
general case. Finally we discuss some examples of Palm
  distributions for specific models and some applications.
\end{abstract}


\section{Introduction}
A spatial point process $\bX$ is briefly speaking a random subset of the
$d$-dimensional Euclidean space $\mathbb R^d$, where $d=2,3$ are the
cases of most practical importance. 
We refer to the (random) elements of $\bX$ as `events' to distinguish
them from other possibly fixed points in $\R^d$.
When studying spatial point
process models and making statistical inference, the conditional
distribution of $\bX$ given a realization of $\bX$ on some specified
region or given the locations of one
or more events in $\bX$ plays an important role, see e.g.\
\cite{moeller:waagepetersen:04} and 
\cite{chiu:stoyan:kendall:mecke:13}. In this paper we focus on the
latter type of conditional distributions which are formally defined in terms of
so-called Palm distributions, first introduced by \citet{palm:43}
for stationary point processes on the real line. Rigorous
definitions and generalizations of Palm distributions to $\mathbb R^d$
and more abstract spaces
have mainly been developed in
probability theory, see \citet{jagers:73} for references and
an historical account. 
Palm distributions are, at least among
many applied
statisticians and among most students, 
considered one of  the more difficult topics in the field
of spatial point processes. This is partly due to the general
definition of Palm distributions which relies on measure theoretical
results,  see e.g.\
\citet{moeller:waagepetersen:04} and
\citet{daley:vere-jones:08} or the references mentioned in Section~\ref{sec:concludingremarks}. The account of conditional distributions
for point processes in \cite{last:90} is mainly
intended for probabilists and is not easily accessible due to an
abstract setting and extensive use of measure theory.

This tutorial provides an introduction to
Palm distributions for spatial point processes. Our setting and
background material on point processes are given in
Section~\ref{s:pre}. Section~\ref{s:finite}, in the context of
  finite point processes, provides an explicit definition of Palm distributions in terms of their density
functions while  Section~\ref{s:general} reviews Palm distributions in the
general case. Section~\ref{s:exs} discusses 
examples of Palm distributions for specific models and
Section~\ref{sec:applications} considers applications of Palm
distributions in the statistical literature.

\section{Prerequisites}\label{s:pre}

\subsection{Setting and notation}

We view a point process as a
random locally finite subset $\mathbf X$ of a Borel set $S \subseteq\mathbb R^d$; 
for measure theoretical details, see e.g.\ 
\citet{moeller:waagepetersen:04} or
\citet{daley:vere-jones:03}.
Denoting $\mathbf X_B=\mathbf X\cap B$
the restriction of $\mathbf X$ to a set $B\subseteq S$,
and $N(B)$ the number of events in $\mathbf X_B$,
local finiteness of $\mathbf X$ means that 
$N(B)<\infty$
almost surely (a.s.) whenever
$B$ is bounded. We denote by $\mathcal B_0$ the family of all bounded
Borel subsets of $S$ and by ${\cal N}$  
the state space consisting of 
the locally
finite subsets (or point configurations) of~$S$. 
Section~\ref{s:finite} considers the case where $S$ is
bounded and hence ${\cal N}$ is all finite subsets of $S$, while 
Section~\ref{s:general} deals with the general case where $S$ is
arbitrary, i.e., including the case $S=\mathbb R^d$.


\subsection{Poisson process}\label{sec:poissonprocess}

The Poisson process is of its own interest and also used for constructing other point
processes as demonstrated in Section~\ref{sec:density} and Section~\ref{s:exs}. 

Suppose $\rho:S\mapsto[0,\infty)$ is
a locally integrable function, that is, $\alpha(B):=\int_B\rho(x)\,\mathrm
dx<\infty$ whenever $B\in\mathcal B_0$. Then   
$\mathbf X$ is a {\it Poisson process} with intensity function $\rho$ if
for any $B\in\mathcal B_0$,
 $N(B)$ is Poisson distributed with mean $\alpha(B)$, 
and conditional on $N(B)=n$, the $n$
events are independent and identically distributed, with a density
proportional to $\rho$ (if
$\alpha(B)=0$, then $N(B)=0$). In fact, this definition is equivalent to
that for any $B\in\mathcal B_0$ and any non-negative measurable function $h$
on $\{ \bx \cap B| \bx \in {\mathcal N}\}$, 
\begin{align}
\E h(\bX_B)=&\,\sum_{n=0}^\infty \frac{\exp\{-\alpha(B)\}}{n!} \nonumber\\
&\,\int_B\cdots\int_B h(\xn)\rho(x_1)\cdots\rho(x_n)\,\mathrm
dx_1\,\cdots\,\mathrm dx_n \,,\label{eq:poissonexp1}
\end{align}
where for $n=0$ the term is
$\exp\{-\alpha(B)\}h(\emptyset)$, where $\emptyset$ is the empty point configuration. 

Note that the definition of a Poisson
process only requires the existence of the intensity measure
$\alpha$, since an event of the process restricted to $B\in\mathcal B_0$
has probability distribution $\alpha(\cdot\cap B)/\alpha(B)$
provided $\alpha(B)>0$. 
We shall use this extension of the definition in Section~\ref{sec:sncpexample}.

\subsection{Finite point processes specified by a density}\label{sec:density}

Let $\bZ$ 
denote a unit rate Poisson process on  $S$, i.e.\ a Poisson process
of constant intensity $\rho(u)=1$, $u \in S$.
Assume that $S$ is bounded and that the 
distribution of $\bX$ is absolutely continuous
with respect to the distribution of $\bZ$ (in short with respect to
$\bZ$) 
with density $f$. Thus, 
for any non-negative measurable function $h$ on ${\cal N}$,
\begin{equation}\label{eq:dens} \E h(\bX) = \E \{f(\bZ) h(\bZ)\}.
\end{equation} 
Moreover, by \eqref{eq:poissonexp1},
\begin{align}
 \E h(\bX) = &\,\sum_{n=0}^\infty \frac{\exp(-|S|)}{n!} \nonumber \\
&\,\int_S\cdots\int_S h(\xn ) f(\xn)\,\mathrm
dx_1\,\cdots\,\mathrm dx_n \label{eq:poissonexp}
\end{align}
where $|S|$ denotes the Lebesgue
measure of $S$. This motivates considering probability statements in terms of
$\exp(-|S|)f(\cdot)$. For example, with $h(\bx)=1(\bx=\emptyset)$, where $1(\cdot)$
denotes the indicator function, we obtain that $P(\bX=\emptyset)$  is
$\exp(-|S|)f(\emptyset)$. Further, for $n \ge 1$,
\[\exp(-|S|)f(\xn) \, \dd x_1\,\cdots\,\mathrm
  dx_n\] 
is  the probability that
  $\mathbf X$ consists of precisely $n$ events with one event
  in each of $n$ infinitesimally small disjoint sets $B_1,\ldots,B_n$ around $\xnl$
with volumes $\mathrm d x_1,\ldots\,\mathrm d x_n$,
respectively. Loosely speaking this is `$P(\bX=\xn)$'.

Suppose we have observed $\bX_ B= \bx_B$ and we wish to predict
the remaining point process $\bX_{S \setminus B}$. Then it is natural
to consider the conditional distribution of $\bX_{S \setminus B}$
given $\bX_B= \bx_B$. By definition of a Poisson process,  $\bZ_B$ and $\bZ_{S\setminus B}$ ($\bZ=\bZ_B
\cup \bZ_{S\setminus B}$) are each
independent unit rate Poisson processes on respectively $B$ and $S
\setminus B$. Thus, in analogy with conditional
densities for multivariate data, this conditional distribution can be
specified in terms of the conditional density
\[ f_{S \setminus B}(\bx_{S \setminus B}|\bx_B) = \frac{ f(\bx_B \cup
  \bx_{S \setminus B})}{f_B(\bx_B)} \]
with respect to $\bZ_{S\setminus B}$ and where 
\[ f_B(\bx_B)= \E f(\bZ_{S \setminus B} \cup \bx_B) \]
is the marginal density of $\bX_B$ with respect to $\bZ_{B}$. Thus the
conditional distribution given a realization of $\bX$ on some prespecified region $B$
is conceptually quite straightforward. Conditioning on that some
prespeficied events belong to $\bX$ is more intricate but an explicit account of this is provided in the next section where it is still assumed that
$\bX$ is specified in terms of  a density.

\section{Palm distributions in the finite case}\label{s:finite}

To understand the definition of a Palm
distribution, it is useful to assume first that $S$ is bounded and that $\bX$ has a density as
introduced in Section~\ref{sec:density} with respect to a unit rate
Poisson process $\bZ$. We make this assumption in the present section, while the general case will be treated in Section~\ref{s:general}.

\subsection{Conditional intensity and joint intensities}

Suppose $f$ is {\it hereditary}, i.e., for any 
pairwise
distinct $x_0,x_1,\ldots,x_n\in S$, 
$f(\xn)>0$ whenever $f(\{x_0,x_1,\ldots,x_n\})>0$. 
We can then define the so-called $n$-th order 
{\it Papangelou conditional intensity} by 
\begin{equation}\label{eq:papangelou}
\lambda^{(n)}(x_1,\ldots,x_n,\bfx)={f(\bfx \cup \xn)}/{f(\bfx)}
\end{equation}
for pairwise distinct $x_1,\ldots,x_n\in S$ and
$\bfx \in\mathcal N\setminus\{x_1,\ldots,x_n\}$, setting
$0/0=0$. By the previous interpretation of $f$,
$\lambda^{(n)}(x_1,\ldots,x_n,\bfx)\dd x_1\cdots\,\mathrm \dd x_n$ can
be considered as the conditional probability of observing
one event in each of the aforementioned infinitesimally small  sets
$B_i$, conditional on that $\bX$ outside $\cup_{i=1}^n B_i$ agrees with
$\bfx$.

For any $n=1,2,\ldots$, we define for pairwise distinct $\xnl\in S$ the $n$-th order {\em joint
intensity function} $\rho^{(n)}$ by
\begin{equation}\label{eq:prodintensdens}
\rho^{(n)}(\xnl)=\mathrm Ef(\mathbf Z\cup \xn) 
\end{equation}
provided the right hand side exists. Particularly, $\rho=\rho^{(1)}$ 
is the usual {\em intensity} function. If $f$ is hereditary, then
$\rho^{(n)}(\xnl)=\mathrm E \lambda^{(n)}(x_1,\ldots,x_n,\bX)$ and by the interpretation of $\ld^{(n)}$ it follows that 
$\rho^{(n)}(\xnl)\,\mathrm
\dd x_1\cdots\,\mathrm \dd x_n$ can be viewed as the probability that $\mathbf X$ has an
event in each of $n$ infinitesimally small sets around $\xnl$
with volumes $\mathrm d x_1,\ldots\,\mathrm d x_n$,
respectively. Loosely speaking, this is `$P(\xnl \in \bX)$'.

Combining \eqref{eq:dens} and \eqref{eq:prodintensdens} 
with either \eqref{eq:poissonexp} or
the extended Slivnyak-Mecke theorem for the Poisson process
given later in \eqref{eq:S-M},
it is straightforwardly seen that
\begin{align}\label{eq:prodintens} 
&\E \sum_{\xnl \in \bX}^{\neq} h(\xnl)\nonumber\\
= &\,\int_{S}\cdots\int_{S}
h(\xnl) \rho^{(n)}(\xnl) \dd x_1 \ldots \dd x_n 
\end{align}
for any non-negative measurable function $h$ on $S^n$, where $\neq$
over the summation sign means that $\xnl$ are pairwise distinct. 
Denoting $N=N(S)$ the number of events in $\bX$, the left hand side in
\eqref{eq:prodintens} with $h=1$ is seen to be the factorial moment 
$\mathrm E\{N(N-1)\cdots(N-n+1)\}$.

\subsection{Definition of Palm distributions in the finite case}

Now, suppose $\xnl\in S$ are
pairwise distinct and $\rho^{(n)}(\xnl)>0$. Then we define the {\em
  reduced Palm distribution} of $\bX$ 
given events at $\xnl$ as the point process 
distribution $\PP^!_{\xnl}$  with
density 
\begin{equation}\label{eq:palmdens} 
f_{\xnl}(\bfx)=\frac{f(\bfx \cup \xn)}
{\rho^{(n)}(\xnl)}, \quad \bfx \in {\cal N},\ \bfx \cap \xn=\emptyset,
\end{equation}
with respect to $\bZ$. We denote by $\bXred$ a point process distributed according to $\PP^!_{\xnl}$.
 If $\xnl\in S$ are not
pairwise distinct or $\rho^{(n)}(\xnl)$ is zero, the choice of $\bXred$ and
its
distribution
$\PP^!_{\xnl}$ is not of any importance for the results in this paper.
Furthermore, the (non-reduced) {\em
  Palm distribution} of $\bX$ 
 given events at $\xnl$ is simply
the distribution of the union $\bXred\cup\xn$. 

\subsection{Remarks}

By the previous infinitesimal
  interpretations of $f$ and $\rho^{(n)}$, 
we can view $\exp(-|S|)f_{\xnl}(\bfx)$  
as the `joint probability' that $\bX$ equals the union  $\bfx \cup \xn$ divided by the
`probability' that $\xnl \in \bX$.
Thus $\PP^!_{\xnl}$ has an
interpretation as the conditional distribution of  $\bX\setminus\xn$
given that $\xnl\in\bX$. 
Conversely, by \eqref{eq:palmdens} with $\bx=\emptyset$ and the remark
just below \eqref{eq:poissonexp},
\begin{equation}\label{eq:likelihoodfact}
\exp(-|S|)f(\{x_1,\ldots,x_n\})=\rho^{(n)}(x_1,\ldots,x_n)
\mathrm P\left(\mathbf X_{\{x_1,\ldots,x_n\}}^! =\emptyset\right)
\end{equation}
provides a factorization into the `probability' of observing
$\{x_1,\ldots,x_n\}$ times the conditional probability of not observing
further events. 

 We obtain immediately from 
\eqref{eq:prodintensdens}
and \eqref{eq:palmdens} that for any pairwise distinct $\xnl\in S$ and $m=1,2,\ldots$, 
$\mathbf X_{\xnl}^!$ has $m$-th order joint intensity function
\begin{equation}\label{eq:palmprodintens}
\rho^{(m)}_{x_1,\ldots,x_n}(u_1,\ldots,u_{m})=\left\{
\begin{array}{ll}
\frac{\rho^{(m+n)}(u_1,\ldots,u_{m},x_1,\ldots,x_n)}{\rho^{(n)}(x_1,\ldots,x_n)} & 
\mbox{if $\rho^{(n)}(x_1,\ldots,x_n)>0$}\\
0 & \mbox{otherwise}
\end{array}
\right.
\end{equation}
for pairwise distinct $u_1,\ldots,u_m \in S\setminus\xn$. 
We write $\rho_{x_1,\dots,x_n}$ for the intensity $\rho^{(1)}_{x_1,\dots,x_n}$.
By \eqref{eq:palmdens}
and \eqref{eq:palmprodintens} we further obtain
\begin{equation}\label{eq:iteratedpalm}
\left(\bX_{x_1,\ldots,x_m}^!\right)_{x_{m+1},\ldots,x_n}^! \stackrel{d}{=} \, \bX_{x_1,\ldots,x_n}^!
\end{equation}
whenever $0<m<n$ and $x_1,\ldots,x_n$ are pairwise distinct, where $\stackrel{d}{=}$ means equality in distribution.

The so-called pair correlation function is for $u,v \in S$ defined
as
\[g(u,v)=\rho^{(2)}(u,v)/\{\rho(u)\rho(v)\}\] 
provided
$\rho(u)\rho(v)>0$ (otherwise we set $g(u,v)=0$). If $\rho(u)\rho(v)>0$, then 
\begin{equation}\label{eq:pcf}
g(u,v)= {\rho_v(u)}/{\rho(u)}={\rho_u(v)}/{\rho(v)},
\end{equation}
cf.\ \eqref{eq:palmprodintens}.
Thus, $g(u,v)>1$
($g(u,v)<1$) means that the presence of an event at $u$ yields an
elevated (decreased) intensity at $v$ and vice versa. 

For later use, notice that 
\begin{align}
  &  \E \sum_{\xnl \in \bX}^{\neq} h(\xnl,\bX
\setminus \xn ) \nonumber \\
= & \int_{S}\cdots\int_{S} \E h(\xnl,\bXred ) \rho^{(n)}(\xnl) \dd x_1
\cdots \dd x_n  \label{eq:palmdef}
\end{align}
for any non-negative measurable function $h$ on $S^n \times {\cal
  N}$. 
 This is called the Campbell-Mecke formula and is straightforwardly verified using \eqref{eq:poissonexp} and
\eqref{eq:palmdens}.
Assuming $f$ is hereditary and rewriting the expectation in the right hand side of 
\eqref{eq:palmdef} in terms of 
\[ f_{\xnl}(\bfx)= 
f(\bfx){\ld^{(n)}(\xnl,\bfx)}/{\rho^{(n)}(\xnl)}\,, \]
the finite point process case of the celebrated 
{\it Georgii-Nguyen-Zessin (GNZ)
formula} 
\begin{align}
  &  \E \sum_{\xnl \in \bX}^{\neq} h(\xnl,\bX
\setminus \xn ) \nonumber \\
= & \int_{S}\cdots\int_{S} \E h(\xnl,\bX ) \lambda^{(n)}(\xnl,\bX)  \dd x_1
\cdots \dd x_n  \label{eq:gnz}
\end{align}
is obtained
 \citep{georgii:76,nguyen:zessin:79b}. We return to
the GNZ formula in connection to Gibbs processes in Section~\ref{sec:gibbsexample}.

\section{Palm distributions in the general case}\label{s:general}

The definitions and results in Section~\ref{s:finite}
extend to
the general case where  $S$ is any Borel subset of $\mathbb
R^d$. 
However, 
 if $|S|=\infty$, the unit rate Poisson process on $S$ will be infinite
and we can not in general assume that $\bX$ is absolutely continuous
with respect to the distribution of this process. Thus we do not
longer have the direct 
definitions \eqref{eq:prodintensdens} and \eqref{eq:palmdens} of
$\rho^{(n)}$ and $\bXred$ in terms of density functions.

\subsection{Definition of Palm distributions in the general case}

Define the $n$-th order factorial moment
measure $\alpha^{(n)}$ on $S^n$ by
\[ \alpha^{(n)}(\times_{i=1}^n B_i ) = \E \sum_{\xnl \in \bX}^{\neq} 1(x_1 \in
B_1, \ldots, x_n \in B_n) \]
for Borel sets $B_1,\ldots,B_n \subseteq S$. Then provided
$\alpha^{(n)}$ is absolutely continuous with respect to Lebesgue
measure on $S^n$, the $n$-th order
joint intensity for $\mathbf X$ is defined as the density of
$\alpha^{(n)}$ with respect to Lebesgue measure on $S^n$. Then, by
standard measure theoretical arguments, \eqref{eq:prodintens} also holds
in the general case. Define further the $n$-th order reduced Campbell
measure by
\[ \CC^!\left( \times_{i=1}^n B_i \times F\right) = \E \sum_{\xnl \in \bX}^{\neq} 1(x_1 \in
B_1, \ldots, x_n \in B_n, \bX \setminus \{\xnl\} \in F) \]
for  Borel sets $B_1,\ldots,B_n \subseteq S$ and any measurable set
$F$ of point configurations in $\cal N$. Obviously for any such $F$,
$C^!$ is dominated by $\alpha^{(n)}$, and so, under suitable
regularity conditions \citep[e.g.\ Section~13.1 in][]{daley:vere-jones:03}, we have
a disintegration of $\CC^!$,
\begin{equation}\label{eq:disintegration} \CC^!( \times_{i=1}^n B_i \times F) = \int_{\times_{i=1}^n B_i}
\PP^!_{\xnl}(F) \alpha^{(n)}(\dd x_1 \cdots \dd x_n ) \end{equation}
where $\PP^!_{\xnl}(F)$ is unique up to an $\alpha^{(n)}$ null-set and
for almost all $\xnl \in S$ defines a distribution of a point process
$\bX^!_{\xnl}$. Again by standard measure theoretical results,
\eqref{eq:palmdef} also holds in the general case (if $\rho^{(n)}$ does
not exist, then  in \eqref{eq:palmdef} replace $\rho^{(n)}(\xnl) \dd x_1 \cdots \dd x_n$ by $\alpha^{(n)}(\dd x_1 \cdots \dd x_n)$). In the finite point
  process case as considered in Section~\ref{s:finite}, by \eqref{eq:palmdef} and \eqref{eq:disintegration} the density approach and the Campbell
  measure approach to define Palm distributions agree.

\subsection{Remarks}

In the general setting, $\rho^{(n)}(\xnl)$ and $\PP^!_{\xnl}$
are clearly only determined up to an $\alpha^{(n)}$
nullset of $S^n$. For simplicity and since there are usually
natural choices of $\rho^{(n)}(\xnl)$ and $\PP^!_{\xnl}$, such nullsets are often
ignored. Further, like in the finite case, $\rho^{(n)}(\xnl)$ and $\PP^!_{\xnl}$ are
  invariant under permutations of the points $x_1,\ldots,x_n$, and
  \eqref{eq:palmprodintens} and
  \eqref{eq:iteratedpalm} also hold in the general case.

Suppose that $\bX$ is {\it stationary}, i.e., its
 distribution is invariant under translations in $\mathbb R^d$ and so
 $S=\mathbb R^d$ (unless $\bX=\emptyset$ which is not a case of
 our interest). This is a specially tractable case, which makes an
 alternative description of Palm distributions possible. Let $\rho$
 denote the constant intensity of $\bX$ and let $o$ denote the
 origin in $\mathbb R^d$. First, we define 
\begin{equation}\label{eq:P0F} \PP^!_o(F)=\frac{1}{\rho|B|}\mathrm
E\sum_{x\in\bX_B}1(\bX \setminus \{x\} -x \in F)
\end{equation}
for any $B \in {\mathcal B}_0$ with $|B|>0$, where by stationarity of $\bX$
  the right hand side does not depend on the choice of
  $B$. Second, we define
\begin{equation}\label{eq:palmstat}
 \PP^!_x(F)= \PP^!_o(F-x) 
\end{equation}
for any $x \in \R^d$. One can then check that the $\PP^!_x$, $x \in \R^d$, defined in this
way satisfy \eqref{eq:disintegration} so that \eqref{eq:palmstat} indeed
defines a Palm distribution, see Appendix~C.2 in
\cite{moeller:waagepetersen:04} for details. Note that \eqref{eq:palmstat} implies that $\bX_x^!-x$ and $\bX_o^!$ are
 identically distributed. The reduced Palm distribution $\PP^!_o$ is
 often interpreted  
 as the `conditional distribution for the further
 events in $\bX$
 given a typical event of $\bX$'. 

\section{Examples of Palm distributions}\label{s:exs}
For some classes of point processes, explicit characterizations of the
Palm distributions are possible. Below we consider Poisson processes, Gibbs processes, log Gaussian Cox
processes (LGCPs), and determinantal point processes which share the property that their Palm distributions of
any order are
again respectively Poisson, Gibbs, LGCPs, and determinantal point processes. We also consider
shot-noise Cox processes, where one point Palm distributions
are not shot-noise Cox processes but have  simple characterizations as cluster
processes. The section is concluded with Tables~\ref{tbl:characteristics1} and~\ref{tbl:characteristics2} which summarize key characteristics for the different model classes.

\subsection{Poisson processes}
In the finite case, by \eqref{eq:poissonexp1}, a Poisson process $\bX$
with intensity function $\rho$ has density 
\[f(\bx) = \exp\left(|S|-\int_S\rho(u)\,\mathrm du\right)\prod_{u
  \in \bx} \rho(u).\] 
  By \eqref{eq:papangelou} and \eqref{eq:prodintensdens} it follows that the $n$-th order Papangelou conditional intensities and the $n$-th order joint intensities agree,
\[ \ld^{(n)}(\xnl,\bx)=\rho^{(n)}(\xnl)=\rho(x_1)\cdots\rho(x_n). \]
Further,  by \eqref{eq:palmdens}, $\bXred\stackrel{d}{=}\bX$.

In the general case, we appeal to the extended Slivnyak-Mecke theorem,
which for a Poisson process $\bX$ with intensity
function $\rho$ states that
\begin{align}
  &  \E \sum_{\xnl \in \bX}^{\neq} h(\xnl,\bX
\setminus \xn ) \nonumber \\
= & \int_{S}\cdots\int_{S} \E h(\xnl,\bX ) \rho(x_1)\cdots\rho(x_n) \dd x_1
\cdots \dd x_n  \label{eq:S-M}
\end{align}
for any non-negative measurable function $h$ on $S^n \times {\cal N}$,
see Theorem~3.3 in 
\cite{moeller:waagepetersen:04} and the references therein.
This implies again that $\rho^{(n)}(\xnl)=\rho(x_1)\cdots\rho(x_n)$
and that $\bX^!_{\xnl}$ is just distributed as $\bX$. In fact, the 
property that $\bX^!_{x}\stackrel{d}{=}\bX$ for all $x\in S$ characterizes the Poisson process, see e.g.\
Proposition~5 in \cite{jagers:73}. Further,
it makes it possible to calculate various useful functional
summaries, see e.g.\ \cite{moeller:waagepetersen:04}, and 
constructions such as stationary Poisson-Voronoi tessellations become
manageable, see \cite{moeller:89a,moeller:94}.

\subsection{Gibbs processes} \label{sec:gibbsexample}
Gibbs processes play an important role in statistical physics and spatial statistics,
see \citet{moeller:waagepetersen:04} and the references therein. Below, for ease of presentation, we consider first a finite Gibbs process.

A finite Gibbs process on a bounded set $S \subset \R^d$ is usually specified in terms of its density
or equivalently
in terms of the Papangelou conditional intensity, where the density is of the form
\[ f(\bx) =
\exp \left \{-\sum_{\by \subseteq \bx} \Phi(\by) \right \} \]
for a so-called potential function $\Phi$ on ${\mathcal N}$, while the Papangelou conditional intensity is
\[
\lambda(u,\bx) = \exp  \left\{-\sum_{\by \subseteq \bx\cup \{u\}:\, u\in \by} \Phi(\by) \right \}.
\]
Here $\exp\{\Phi(\emptyset)\}$ is the normalizing constant (partition
function) of $f(\cdot)$ which in general is not expressible on closed form, while $\lambda(u,\bx)$ does not depend on the normalizing constant.
It follows
that the $n$-th order Palm distribution of a Gibbs process with respect to
$\xnl$ is itself a Gibbs
process with potential function $\Phi_{\xnl}(\by)= \Phi(\xn \cup
\by)$ for $\by\not=\emptyset$. Moreover, for pairwise distinct $u_1,\ldots,u_m,x_1,\ldots,x_n \in S$ and $\bx
\in \mathcal N\setminus  \{u_1,\ldots,u_m,x_1,\ldots,x_n\}$, the $m$-th order Papangelou conditional intensity of $\bXred$
is simply 
\[ \lambda^{!(m)}_{\xnl}(u_1,\ldots,u_m,\bx) = \lambda^{(m)}(u_1,\ldots,u_m,\bx \cup \xn) .\]
 

For instance, a first order inhomogeneous pairwise interaction Gibbs point process has first order potential $\Phi(\{u\})=\Phi_1(u)$, second order potential $\Phi(\{u,v\})=\Phi_2(v-u)$, and $\Phi(\by)=0$ whenever the cardinality of $\by$ is larger than
two; see \citet{moeller:waagepetersen:04} for
conditions on the functions $\Phi_1$ and $\Phi_2$ ensuring that the model is
well-defined. The Strauss model \citep{strauss:75,kelly:ripley:76} is
a particular case with $\Phi_1(u)=\theta_1\in \mathbb R$ and
$\Phi_2(u-v)= \theta_2 1(\|u-v\| \leq R)$, for $\theta_2 \geq 0$ and
$0<R<\infty$. The Palm process $\bXred$ becomes
again an inhomogeneous pairwise interaction Gibbs process with inhomogeneous first order potential
$\Phi_{x_1,\ldots,x_n}(\{u\})=\Phi_1(u)+ \sum_{i=1}^n\Phi_2(u-x_i)$ and second order potential
identical  to that of $\bX$. 



In the general case, a Gibbs process can be defined
\citep{nguyen:zessin:79b} in  terms of the GNZ
formula \eqref{eq:gnz} briefly discussed at
the end of Section~\ref{s:finite}: $\bX$ is a Gibbs point process 
with Papangelou conditional intensity $\lambda$
if $\lambda$ is a non-negative measurable function 
on $S\times\mathcal N$ such that
\begin{equation}\label{e:Gibbs}
\E \sum_{x \in \bX} h(x,\bX\setminus \{x\})
=  \E \int_{S}\lambda(x,\bX) h(x,\bX)\dd x
\end{equation} 
for any non-negative measurable function $h$ on $S\times\mathcal
N$. For conditions
ensuring that \eqref{e:Gibbs} {holds, we refer to \cite{ruelle:69}, \cite{georgii:88}, or \citet{dereudre:drouilhet:georgii:12}.

By the extensions of \eqref{eq:prodintens} and
\eqref{eq:palmdef} to the general case, \eqref{e:Gibbs}
implies $\rho(x)=\E\lambda(x,\bX)$. Unfortunately, in general it is not
feasible to express $\rho(x)=\E\lambda(x,\bX)$ on closed form, though
approximations exist \citep{baddeley:nair:12}. 
Also, for Gibbs processes, the pair correlation function $g(u,v)$ can be below or above 1 depending on $u$ and $v$ \citep[see e.g.\ pages 240-241
  in][]{illian:et:al:08}, and so from~\eqref{eq:pcf}, $\rho_v(u)$ may
be smaller or larger than $\rho(u)$, depending on $u$ and $v$.}
Moreover, for pairwise distinct $\xnl\in S$,
$\mathrm P^!_{\xnl}$ is
absolutely continuous with respect to the distribution of $\bX$, with
density $\tilde f(\bx) = \lambda^{(n)}(\xnl,\bx)/\rho^{(n)}(\xnl)$, where
\begin{align*}
\lambda^{(n)}(\xnl,\mathbf x)=&\,\lambda(x_1,\mathbf x)
\lambda(x_2,\mathbf x\cup\{x_1\})\\
&\cdots
\lambda(x_n,\mathbf x\cup\{x_1,\ldots,x_{n-1}\}) 
\end{align*}
for $\xnl\in S$ and $\mathbf x\in\mathcal N$. This follows from
\eqref{eq:gnz} and \eqref{e:Gibbs} and is in accordance with
\eqref{eq:papangelou} and \eqref{eq:iteratedpalm}. Note that in
  this connection, the roles of $\bx$ and $\xnl$ in
  $\lambda^{(n)}(\xnl,\bx)$ are interchanged: now
  $\xnl$ are fixed while $\bx$ is the variable argument of the density
  of $\PP^!_{\xnl}$.

\subsection{Cox processes}\label{sec:cox}
Let $\mathbf \Lambda=\{\Lambda(x)\}_{x\in S}$ be a 
non-negative random field such that  $\mathbf \Lambda$ is locally
integrable a.s., that is,
for any $B\in\mathcal B_0$, the integral
$\int_B \Lambda(x)\,\mathrm dx$ exists and is finite a.s. Suppose
$\mathbf X$ is a Cox process with
random intensity function $\mathbf \Lambda$, i.e.,
conditional on $\mathbf
\Lambda$, $\mathbf X$ is a Poisson process with intensity function
$\mathbf \Lambda$. Apart from very simple models of $\mathbf \Lambda$ 
such as all $\Lambda(x)$ being equal to the same random
variable following e.g.\ a gamma distribution, 
the density of $\mathbf X$ restricted to
a set $B\in\mathcal B_0$ is intractable. However, if
$\mathbf \Lambda$ has moments
of any order $n=1,2,\ldots$, then by conditioning on $\mathbf\Lambda$ 
we immediately obtain
\begin{equation}\label{eq:coxprodintens}
\rho^{(n)}(x_1,\ldots,x_n)=\mathrm
E\left\{\prod_{i=1}^n\Lambda(x_i)\right\}
\end{equation} 
for any pairwise distinct 
$x_1,\ldots,x_n\in S$. For any  $B\in\mathcal B_0$ the conditional
density of $\bX \cap B$ given $\mathbf\Lambda$ is
\[ f(\bx | \mathbf\Lambda) = \exp\left (|B|-\int_B \Lambda(u) \dd u \right ) \prod_{u \in \bx} \Lambda (u) \]
and it follows that the marginal density and the reduced Palm density
of $\bX \cap B$ are given by
\[ f(\bx)= \E f(\bx | \mathbf\Lambda) \quad \text{ and } \quad f_\xnl (\bx)= \E
\left \{ f(\bx | \mathbf\Lambda) \frac{\prod_{i=1}^n \Lambda(x_i)}{
    \rho^{(n)}(\xnl)} \right \}.\]
The expression for the reduced Palm density in fact shows that the
reduced Palm distribution of $\bX \cap B$ is also a Cox process but now with
a random intensity function $\mathbf\Lambda_\xnl$ that has density
\[ \frac{\prod_{i=1}^n \Lambda(x_i)}{ \rho^{(n)}(\xnl)} \]
with respect to the distribution of $\mathbf\Lambda$, i.e.\
\[ P(\mathbf{\Lambda}_\xnl \in A) = \E \left \{  1(\mathbf\Lambda \in A) \frac{\prod_{i=1}^n \Lambda(x_i)}
{\rho^{(n)}(\xnl)} \right \}\]
for subsets $A$ of the sample space of $\mathbf \Lambda$. The density
perspective gives a very simple derivation of this result which in
fact also holds for general Cox processes, see e.g.\ Example~13.1(a) in \cite{daley:vere-jones:08} or page 169 in \cite{chiu:stoyan:kendall:mecke:13}.

More generally, conditioning on $\mathbf\Lambda$ and
using \eqref{eq:palmdef} and \eqref{eq:S-M},
the reduced Palm distributions satisfy
\begin{align}
\E&\left\{h\left(\xnl,\bXred\right)\right\}  \rho^{(n)}(x_1,\ldots,x_n)\nonumber\\
&=\E\left\{h(\xnl,\bX)\prod_{i=1}^n\Lambda(x_i)\right\} 
\label{e:cox1}
\end{align}
for a non-negative measurable function $h$ on $S^n \times \mathcal N$.
In the sequel, we consider
distributions of $\mathbf \Lambda$, where
\eqref{eq:coxprodintens}-\eqref{e:cox1} become useful. 

\subsubsection{Log Gaussian Cox  processes}\label{sec:lgcpexample}

Let $\Lambda(x)=\exp\{Y(x)\}$, where 
$\mathbf Y=\{Y(x)\}_{x\in S}$ is a Gaussian process with
mean function $\mu$ and covariance function $c$ 
so that 
$\mathbf\Lambda$ is locally integrable a.s. 
\citep[simple conditions ensuring this 
are given in][]{moeller:syversveen:waagepetersen:98}.
 Then $\mathbf X$ is a log Gaussian Cox
process (LGCP) as introduced by \citet{coles:jones:91} in astronomy  and
independently by \citet{moeller:syversveen:waagepetersen:98} in
statistics.
By \citet[][Theorem~1]{moeller:syversveen:waagepetersen:98},
for pairwise distinct $x_1,\ldots,x_n\in S$,
\begin{equation}\label{eq:prodintenslgcp}
\rho^{(n)}(x_1,\ldots,x_n)=\left\{\prod_{i=1}^n\rho(x_i)\right\}
\left\{\prod_{1\le i<j\le n}g(x_i,x_j)\right\},
\end{equation}
where $\rho(x)=\exp\{\mu(x)+c(x,x)/2\}$ is the intensity function and
the pair correlation function \eqref{eq:pcf} is 
$g(u,v)=\exp\{c(u,v)\}$.  The intensity function of $\bXred$ takes the form
\begin{equation}\label{eq:lgcpintens}
  \rho_{\xnl}(u) = \rho(u) \prod_{i=1}^n g(u,x_i)
\end{equation}
so in the common case where $c$ is positive, the intensity of
$\bXred$ is larger than that of $\bX$.

In \cite{coeurjolly:moeller:waagepetersen:15} it is verified that 
for pairwise distinct \linebreak$\xnl\in S$, 
$\bXred$ is an LGCP with underlying Gaussian process
$\{Y(x)+\sum_{i=1}^nc(x,x_i)\}_{x\in S}$. Note that this 
 Gaussian process
also has covariance function $c$ but its mean
  function is $\mu_{\xnl}(x)=\mu(x)+\sum_{i=1}^nc(x,x_i)$. 
\cite{coeurjolly:moeller:waagepetersen:15} discuss how this result can be exploited for
functional summaries. Moreover, if the covariance function $c$ is
non-negative, $\bX$ is distributed as an independent thinning of
$\bXred$ with inclusion probabilities
$p(x)=\exp \{-\sum_{i=1}^nc(x,x_i)\}$.

\subsubsection{Shot noise Cox processes}\label{sec:sncpexample} 
For a shot noise Cox process \citep{moeller:02},
\[\Lambda(x)=\sum_j\gamma_j k(c_j,x),\]
where $k(c_j,\cdot)$ is a kernel (i.e., a density function for a
continuous $d$-dimensional random variable) and the $(c_j,\gamma_j)$
are the events of a Poisson process $\mathbf\Phi$ on $\mathbb
R^d\times(0,\infty)$ with intensity measure $\alpha$ 
so that $\mathbf\Lambda$ becomes locally integrable
a.s. 
It can be viewed as a cluster process $\bX=\cup_{j} \bY_j$, where
conditional on $\mathbf\Phi$, the cluster $\bY_j$ is a Poisson process
with intensity function $\gamma_j k(c_j,\cdot)$ and the clusters are
independent.

The intensity function
is
\[\rho(x)=\int \gamma k(c,x)\,\mathrm d\alpha(c,\gamma),\]
provided the integral is finite for all $x\in S$. Making this
assumption, it can be verified \citep[Proposition 2 in][]{moeller:02} that for $x\in S$ with $\rho(x)>0$,
$\bX^!_x$ is a Cox process with random intensity
function $\Lambda(\cdot)+\Lambda_x(\cdot)$, where 
$\Lambda_x(\cdot)=\gamma_x k(c_x,\cdot)$, and where $(c_x,\gamma_x)$ is a random
variable independent of $\mathbf\Phi$ and 
defined on $S\times(0,\infty)$ such that for any Borel set
$B\subseteq S\times(0,\infty)$, 
\[\mathrm P\left\{(c_x,\gamma_x)\in B\right\}=\frac{\int_B\gamma k(c,x)
\,\mathrm d\alpha(c,\gamma)}{\rho(x)}.\]
In other words, $\bX^!_x$ is
distributed as $\bX\cup\bY_x$, where $\bY_x$ is independent of $\bX$
and conditional on $(c_x,\gamma_x)$, the `extra cluster' 
$\bY_x$ is a finite Poisson
process with intensity function $\gamma_xk(c_x,\cdot)$. Thus, like for an
LGCP with positive covariance function, 
$\bX^!_x$ has a higher
intensity than $\bX$.

For instance, if 
$\mathrm d\alpha(c,\gamma)=\mathrm dc\,\mathrm d\chi(\gamma)$, where
$\chi$ is a locally finite measure on $(0,\infty)$, then 
$\rho(x)=\kappa f(x)$, where it is
assumed that 
$\kappa=\int\gamma\,\mathrm d\chi(\gamma)<\infty$ and
$f(x)=\int k(c,x)\,\mathrm dc<\infty$, and furthermore, for $\rho(x)>0$,
$c_x$ and
$\gamma_x$ are independent, $c_x$ follows the density 
$k(\cdot,x)/f(x)$, and $\mathrm P(\gamma_x\in A)=\kappa^{-1}
\int_A\gamma\,\mathrm d\chi(\gamma)$. The special case
of a Neyman-Scott process \citep{neyman:scott:58} occurs when
$S=\mathbb R^d$, 
$\chi$ is concentrated at a given value $\gamma>0$, $\chi(\{\gamma\})<\infty$, and
$k(c,\cdot)=k_o(\cdot-c)$, where $k_o$ is a density function. Then 
$\bX$ is stationary, $\rho=\kappa=\gamma\chi(\{\gamma\})$, $c_x$ has
density $k_o(x-\cdot)$, and conditional on $c_x$, 
$\bY_x$ is a finite Poisson
process with intensity function $\gamma k_o(\cdot-c_x)$. 
Examples include a (modified) Thomas
process, where $k_o$ is a zero-mean normal density, and a Mat{\'e}rn
cluster process, where $k_o$ is a uniform density on a ball centered
at the origin. For $n>1$, 
the $n$-th order reduced Palm distributions become more complicated. 

In a Neyman-Scott process, the number of events in the clusters are
independent and identically Poisson distributed. For a general stationary Poisson
cluster process the cluster centres still form a stationary Poisson process but the Poisson distribution of the number of events in a cluster is replaced by
any discrete distribution on the non-negative integers. Finally, we
notice that the Palm
distribution for stationary Poisson
cluster processes and more generally infinitely divisible point
processes can also be derived, see 
\citet{chiu:stoyan:kendall:mecke:13} and the references therein.

\subsection{Determinantal point processes} \label{sec:dpp}
Determinantal point processes is a class of repulsive point processes
that has recently attracted interest for statistical applications, see
\cite{lavancier:moeller:rubak:15} and the references therein. For simplicity we restrict attention to determinantal point processes specified by a covariance function
$C:\mathbb R^d\times\mathbb R^d\mapsto\mathbb C$ such that $\int_S C(u,u)\,\mathrm d u<\infty$ whenever $S\subset\mathbb R^d$ is compact. Then $\bX$ is said to be a
 determinantal point process with kernel $C$ if for any $n=1,2,\ldots$ and pairwise distinct $x_1,\ldots,x_n\in\mathbb R^d$, the $n$-th order joint intensity function exists and is given by
\begin{equation}\label{e:DPPdef}
\rho^{(n)}(x_1,\ldots,x_n)={\mathrm{det}}[C](x_1,\ldots,x_n)
\end{equation}
where $[C](x_1,\ldots,x_n)$ denotes the matrix with entries $C(x_i,x_j)$, $i,j=1,\ldots,n$. The determinant of this matrix will not depend on the ordering of $x_1,\ldots,x_n$, so we also write ${\mathrm{det}}[C](\{x_1,\ldots,x_n\})$ for ${\mathrm{det}}[C](x_1,\ldots,x_n)$. Note that $\rho(u)=C(u,u)$ is the intensity function.  

The existence of the process is equivalent to that for any compact set $S\subset\mathbb R^d$,
the eigenvalues of the kernel restricted to $S\times S$ are at most 1. The process is then uniquely characterized by \eqref{e:DPPdef}. If the eigenvalues are strictly less than 1, then $\bX$ restricted to $S$ has density
\[f(\bx)\propto{\mathrm{det}}[\tilde C](\bx) \]
where $\tilde C$ is the covariance function given by the integral equation
\[\tilde C(x,y)-\int_S\tilde C(x,z)C(z,y)\,\mathrm dz=C(x,y),\qquad x,y\in S,\] 
and where the normalizing constant of the density can be expressed in terms of the eigenvalues. For further details, see \cite{lavancier:moeller:rubak:15}. 

Consider any pairwise distinct $x,u_1,\ldots,u_n\in\mathbb R^d$ with $\rho(x)>0$, and define the covariance function $C_x$ by 
\[C_x(u_1,u_2)=C(u_1,u_2)-C(u_1,x)\tilde C(x,u_2)/\tilde C(u,u).\] 
Using \eqref{eq:palmprodintens} it follows that $\bX^!_x$ has $n$-th order joint intensity function
\[\rho_x(u_1,\ldots,u_n)={\mathrm{det}}[C_x](u_1,\ldots,u_n).\]
Consequently $\bX^!_x$ is a determinantal point process with kernel $C_x$. 
See also Theorem~6.5 in
\cite{shirai:takahashi:03} or Appendix~C of the supplementary material
for \cite{lavancier:moeller:rubak:15}. By \eqref{eq:iteratedpalm} and induction
it follows that determinantal point processes are closed under Palm
conditioning: the reduced Palm distribution of any order of a
determinantal point process is again a determinantal point process.

\begin{center}
  \begin{table}[htbp]\centering
\begin{tabular}{lp{2.5cm}p{5.5cm}}
\hline 
 Characteristic & Poisson &  Gibbs \\
 \hline
 &&\\
 Density & $z^{-1}_S \prod_{v\in \bx} \rho(v)$ & $\propto\exp \left \{-\sum_{\emptyset\not=\by \subseteq \bx} \Phi(\by) \right \}$\\
 $f(\bx)$ &&\\
 &&\\
 Papangelou cond. &$\rho(u)$& $\exp  \left\{-\sum_{\by \subseteq \bx\cup \{u\}: u\in \by} \Phi(\by) \right \}$\\
 intensity $\lambda(u,\bx)$ &&\\
 &&\\
 Joint intensity &$ \prod_{i=1}^n \rho(x_i)$& $\E f(\{x_1,\dots,x_n\} \cup \bZ)$\\
 $\rho^{(n)}(x_1,\dots,x_n)$ &&\\
 &&\\
 One-point Palm &$z^{-1}_S \prod_{v\in \bx} \rho(v)$& $\propto\exp  \left\{-\sum_{\emptyset\not=\by \subseteq \bx\cup \{u\}} \Phi(\by) \right \}$\\
 density $f_u(\bx)$ &&\\
&&\\
One-point Palm  &$\rho(u)$&$\E f(\{u,v\}\cup\bZ)/\E f(\{v\}\cup\bZ) $\\
intensity $\rho_v(u)$&&\\
&&\\\hline
\end{tabular}  
\caption{Point process characteristics for Poisson and Gibbs processes
when the state space $S$ is bounded.
For the Poisson process, $\rho(\cdot)$
  denotes the intensity function and $z_S=\exp(|S| - \int_S \rho(v) \dd v)$ the normalizing constant. For the Gibbs process,
  $\Phi$ denotes the potential function, $\bZ$ is the unit rate Poisson process on $S$, and the normalizing constants of the density and one-point Palm density and the expectations for the $n$-th order joint intensity and the one-point Palm intensity are in general intractable, see Section~\ref{sec:gibbsexample} for more details.}\label{tbl:characteristics1}
\end{table} 
\end{center}

\begin{center}
  \begin{table}[htbp]\centering
\begin{tabular}{lp{4.7cm}p{4.5cm}}
\hline 
 Characteristic & Cox &  Determinantal \\
 \hline
 &&\\
 Density & $\E f(\bx\mid\mathbf \Lambda)$ &$\propto \det [\tilde C](\bx)$\\
 $f(\bx)$ &&\\
 &&\\
 Papangelou cond. &${\displaystyle {\E f(\bx \cup\{u\}\mid \mathbf \Lambda)}/{\E f(\bx\mid\mathbf \Lambda)}}$ &
 ${\displaystyle {\det [\tilde C](\bx\cup\{u\})}/{\det [\tilde C](\bx)}
 }$ 
 \\
 intensity $\lambda(u,\bx)$ &&\\
 &&\\
 Joint intensity & $\E \prod_{v\in \bx}  \Lambda(v)$& $\det [ C](x_1,\ldots,x_n)$\\
 $\rho^{(n)}(x_1,\dots,x_n)$ &&\\
 &&\\
 One-point Palm &$\E f(\bx \mid \mathbf \Lambda_u)$& $\propto \det [\tilde C_u](\bx)$\\
 density $f_u(\bx)$ &&\\
&&\\
One-point Palm &&\\
intensity $\rho_v(u)$ &$\E\{ \Lambda(u)  \Lambda(v) \}/\E
 \Lambda(v)$&$\det[C](u,v)/C(v,v)$\\
&&\\\hline
\end{tabular}  
\caption{Point process characteristics for Cox and determinantal point processes 
when the state space $S$ is 
compact.
For the Cox process, $\mathbf\Lambda$
  denotes the random intensity function, $\mathbf \Lambda_u$ is the modified random field (see Section~\ref{sec:cox}), and $f(\cdot|\mathbf\Lambda)$ is a Poisson process density when we condition on that $\mathbf\Lambda$ is the intensity function; 
all the expectations are in general intractable.  
  For the determinantal point process, $C$ denotes its kernel and we refer to Section~\ref{sec:dpp} for details on the related kernels $\tilde C$ and $\tilde C_u$; the normalizing constants of the densities are known (see Section~\ref{sec:dpp}).}\label{tbl:characteristics2}
\end{table}
\end{center}

\section{Examples of applications}\label{sec:applications}

In this section we review a number of applications of Palm
distributions in spatial statistics.

\subsection{Functional summary statistics}

Below we briefly consider two popular functional summary statistics, which are used for exploratory purposes as well as model fitting and model assessment.

First, suppose $\bX$ is stationary, with intensity $\rho>0$.
The nearest-neighbour distribution function 
$G$ is defined by $G(t)=\PP^!_o \{\bX \cap b(o,t)
\neq \emptyset\}$, where $b(o,t)$ is the ball centered at $o$ and of
radius $t>0$. Thus $G(t)$ is interpreted as the probability of having an event within distance $t$ from a typical event. Moreover, 
Ripley's $K$-function \citep{ripley:76} times $\rho$ is defined by 
$\rho K(t) = \E \sum_{v \in \bX^!_o} 1(\|v\| \le t)$, that is,
the expected number of further events within distance
$t$ of a typical event.

Second, if the pair correlation function $g(u,v)=g_0(v-u)$ only depends on
$v-u$ (see \eqref{eq:pcf}), the definition of the $K$-function can be extended: The
inhomogeneous $K$-function
\citep{baddeley:moeller:waagepetersen:00} is defined by
\[ K(t)=\int_{\|v\|\le t} g_0(v) \dd v. \]
By \eqref{eq:pcf}, it follows that
\[ K(t) = 
\E \sum_{v \in \bX^!_u} \frac{1(\|v-u\| \le t)}{\rho(v)} \]
for any $u \in S$ with $\rho(u)>0$. If for $\|v-u\| \le
t$, $\rho(v)$ is close to $\rho(u)$,  we obtain $\rho(u) K(t) \approx \E \sum_{v \in \bX^!_u} 1(\|v-u\|
\le t)$. This is a `local' version of the interpretation of $K(t)$ in the
stationary case. 

 Nonparametric estimation of $K$ and $G$ is based on empirical
 versions obtained from \eqref{eq:P0F}. For some parametric Poisson
 and Cox process models, $K$ or $G$ are expressible on closed form and may be compared with corresponding nonparametric estimates when finding parameter estimates or assessing a fitted model. 
 See \citet{moeller:waagepetersen:07} and the references therein. 


\subsection{Prediction given partial observation of a point process}

Suppose $S$ is bounded and we observe a point process $\bY$ contained in 
a finite point process $\bX$
specified by some density $f$ with respect to the unit rate Poisson
process $\bZ$. If $B \subset S$ with $|B|>0$ and $\bY= \bX_B$, then prediction of 
$\bX_{S \setminus B}$ given $\bY=\by$ can be based on the conditional density $f_{S
  \setminus B}(\cdot|\by)$ introduced in
Section~\ref{sec:density}. On the other hand, if we just know that $\by
\subseteq \bX$, then it could be tempting to try to predict $\bX
\setminus \by$ using $\bX^!_\by$. This would in general be
incorrect. For instance, for an LGCP  with positive
covariance function, the intensity of $\bX^!_\by$ can be much larger
than the one of $\bX$, cf.\ \eqref{eq:lgcpintens}. Thus on average $\bX^!_\by \cup \by$ would
contain more events than $\bX$. The issue here is that the reduced Palm
distribution is concerned with the conditional distribution of $\bX$
conditional on that {\em prespecified} points fall in $\bX$. Hence the
sampling mechanism that leads from $\bX$ to $\bY$ must be taken into
account. For instance, if the distribution of $\bY$ conditional on $\bX=\bx$ is
specified by a probability density function $p(\cdot|\bx)$ (on the set of all subsets of $\bx$), then by Proposition~1 in
\cite{baddeley:moeller:waagepetersen:00}, the marginal  density of $\bY$
with respect to $\bZ$ is 
\[ g(\by)= \rho^{(n)}(\by)\exp(|S|) \E \left\{p(\by|\bX_\by^!
  \cup \by) \right\},  \]
where $n=n(\by)$ is the cardinality of $\by$. Thus the conditional
distribution of $\bX \setminus \by$ given $\bY= \by$ has density
\[ f(\bx|\by) = p(\by|\bx \cup \by)f(\bx \cup \by) \exp(|S|)/ g(\by) \]
with respect to $\bZ$. 

\subsection{Mat{\'e}rn-thinned Cox processes}

Some applications of spatial point processes require models that
combine clustering at a large scale with regularity at a local
scale \citep{lavancier:moeller:15}. \cite{andersen:hahn:15} study a class of so-called Mat{\' e}rn
thinned Cox processes where (clustered) Cox processes are subjected to
dependent Mat{\' e}rn type II thinning \citep{matern:86} that introduces regularity in
the resulting point processes. The intensity function and second-order joint
intensity of the Mat{\' e}rn-thinned Cox process is expressed in terms
of univariate and bivariate inclusion probabilities which in turn are
expressed in terms of one- and two-point Palm probabilities for the underlying Cox process extended with a uniformly distributed mark for each event. In case of
an underlying shot-noise Cox process, explicit expressions for the
univariate inclusion probabilities are obtained using the simple
characterization of one-point Palm distributions described in Section~\ref{sec:sncpexample}.

\subsection{Palm likelihood}

Minimum contrast estimators based on the $K$-function or the pair
correlation function or composite likelihood methods are standard
methods to fit parametric models \citep[see e.g.][]{jolivet:91,guan:06,moeller:waagepetersen:07,waagepetersen:guan:09,biscio:lavancier:15}. 
\citet{tanaka:ogata:stoyan:08} proposed an approach
based on Palm intensities to fit parametric stationary
models, which is briefly presented below. 

Given a parametric model $g(u,v)=g_0(v-u;\theta)$ for the pair
correlation function of $\bX$ and a location $u \in S$, the intensity
function of $\bX^!_u$ is $\rho_u(v;\ta)=\rho g_0(v-u;\theta)$
where $\rho$ is the constant intensity of $\bX$ assumed here to be known.
Following \cite{schoenberg:05}, the so-called log composite likelihood score
\[ \sum_{v \in \bX_u^! \cap b(u,R)} \frac{\dd}{\dd \ta} \log \rho_u(v;\ta) -
  \int_{b(u,R)} \rho_u(v;\ta) \dd v \]
forms an unbiased estimating function for $\ta$, where $R>0$ is a user-specified
tuning parameter. Usually $\bX_u^!$ is not known. However, suppose that $\bX$ is observed on $W \in
{\mathcal B}_0$ and in order to introduce a border correction let $W \ominus R= \{ u \in W | b(u,R) \subseteq W\}$. Then, by \eqref{eq:P0F},
\begin{equation}\label{eq:palmscore} \sum_{\substack{u \in \bX \cap W
      \ominus R,\\ v \in \bX \cap b(u,R)}}^{\neq}
\frac{\dd}{\dd \ta} \log \rho_u(v;\ta) - N (W \ominus R) \int_{b(o,R)}
\rho_o(v;\ta) \dd v 
\end{equation}
is an unbiased estimate of the above composite likelihood
score times $\rho |W \ominus R|$. \citet{tanaka:ogata:stoyan:08} coined the antiderivative of
\eqref{eq:palmscore} the Palm likelihood. Asymptotic properties of
Palm likelihood parameter estimates are studied by
\cite{prokesova:jensen:13} who also proposed the border correction
applied in \eqref{eq:palmscore}.


\section{Concluding remarks}\label{sec:concludingremarks}
  
The intention of this paper was to give a brief and non-technical
introdution to Palm distributions for spatial point processes. For
more extensive treatments of the topic we refer to \citet{cressie:93},
\citet{baddeley:99}, \Citet{lieshout:00}, \citet{daley:vere-jones:08},
\citet{chiu:stoyan:kendall:mecke:13}, and \citet{spodarev:13}. We omitted the case of marked point processes for sake of brevity. The theory of Palm distributions for marked point processes is fairly similar to that for ordinary point processes. Accounts of Palm distributions for marked point processes can be found in \citet{chiu:stoyan:kendall:mecke:13} and \citet{heinrich:13}, while summary statistics related to Palm distributions for marked point processes are reviewed in \citet{illian:et:al:08} and \citet{baddeley:10}.

 We finally note that consideration of space-time point processes \citep{diggle:gabriel:10} suggests yet another useful notion of conditioning on the past. For a space-time point process, the conditional intensity for a time-space point $(t,x)$ usually refers to the conditional probability of observing an event at spatial location $x$ at time $t$ given the history of the space-time point process up to but not including time $t$. So this conditional intensity naturally takes the time-ordering into account, while there is no natural ordering when considering a spatial point process.

\subsubsection*{Acknowledgments}
We thank the editor and two referees for detailed and
constructive comments that helped to improve the paper. J. M{\o}ller and R. Waagepetersen are supported by the Danish Council for Independent Research | Natural
Sciences,
grant 12-124675,
"Mathematical and Statistical Analysis of Spatial Data", and by
the "Centre for Stochastic Geometry and Advanced Bioimaging",
funded by grant 8721 from the Villum Foundation. J.-F. Coeurjolly is supported by ANR-11-LABX-0025 PERSYVAL-Lab (2011, project OculoNimbus).

\bibliographystyle{royal} 
\bibliography{palm}


\end{document}